\documentclass[11pt]{article}
\usepackage{geometry}                
\usepackage{amsmath,amssymb}
\usepackage{amsthm}
\usepackage{hyperref}
\usepackage{graphicx,epsfig,color}
\usepackage{booktabs,bm,multirow}
\usepackage{cases}

\newtheorem{theorem}{Theorem}

\newtheorem{remark}[theorem]{Remark}
\newtheorem{lemma}[theorem]{Lemma}

\input{siam.sty}

\title{A doubly stochastic block Gauss--Seidel algorithm for solving linear equations}
\author{Kui Du\thanks{School of Mathematical Sciences, Xiamen University, Xiamen 361005, China ({\tt kuidu@xmu.edu.cn}).},\quad Xiao-Hui Sun\thanks{School of Mathematical Sciences, Xiamen University, Xiamen 361005, China ({\tt 19020190154621@stu.xmu.edu.cn}).}} 
\date{}

\begin{document}
\maketitle

\begin{abstract}  
\vspace{.5mm} 

We propose a simple doubly stochastic block Gauss--Seidel algorithm for solving linear systems of equations. By varying the row partition parameter and the column partition parameter of the coefficient matrix, we recover the Landweber algorithm, the randomized Kaczmarz algorithm, the randomized Gauss--Seidel algorithm, and the doubly stochastic Gauss--Seidel algorithm. For general (consistent or inconsistent) linear systems, we show the exponential convergence of the {\it norms of the expected iterates} via exact formulas. For consistent linear systems, we prove the exponential convergence of the {\it expected norms of the error and the residual}. Numerical experiments are given to illustrate the efficiency of the proposed algorithm.


\end{abstract}

\section{Introduction}
\label{intro}
Randomized iterative algorithms for solving a linear system of equations \begin{equation}\label{lineq}{\bf Ax=b},\quad \mbf A\in\mbbr^{m\times n}, \quad \mbf b\in\mbbr^m\end{equation} have attracted much attention recently; see, for example, \cite{strohmer2009rando,leventhal2010rando,zouzias2013rando,needell2014paved,needell2015rando,ma2015conve,gower2015rando,ma2018itera,bai2018conve,bai2018greed,bai2018relax,bai2019parti,bai2019greed,du2019tight,razaviyayn2019linea,necoara2019faste,zhang2019new,liu2019varia,wu2020proje,chen2020error,niu2020greed,richtarik2017stoch}.
At each step, to generate the next iterate from the current iterate, the randomized Kaczmarz algorithm \cite{strohmer2009rando} uses a randomly picked row, the randomized Gauss--Seidel (i.e., randomized coordinate descent) algorithm \cite{leventhal2010rando} uses a randomly picked column, and the doubly stochastic Gauss--Seidel algorithm \cite{razaviyayn2019linea} uses a randomly picked entry of the coefficient matrix $\mbf A$. It is natural to ask whether one can design a randomized algorithm which uses a randomly picked submatrix of $\mbf A$.

In this paper, we propose a doubly stochastic block Gauss--Seidel (DSBGS) algorithm which uses a submatrix of $\mbf A$ at each step (see Algorithm 1 in \S 2). We can view DSBGS as a stochastic gradient descent for solving the following optimization problem \begin{equation}\label{opt}\min_{\mbf x\in\mbbr^n}\l\{f(\mbf x):=\frac{1}{2\|\mbf A\|_\rmf^2}\|\mbf b-\mbf A\mbf x\|_2^2\r\}.\end{equation} The Landweber iterative algorithm \cite{landweber1951itera}, the randomized Kaczmarz (RK) algorithm, the randomized Gauss--Seidel (RGS) algorithm, and the doubly stochastic Gauss--Seidel (DSGS) algorithm are special cases of our algorithm. Our algorithm does not need to use projections and Moore-Penrose pseudoinverses of submatrices, so it is different from the block algorithms in \cite{needell2014paved,needell2015rando,gower2015rando}. Numerical experiments for both synthetic data and real-world data are given to illustrate the efficiency of DSBGS.  

{\it Main contributions}. We propose a simple doubly stochastic block Gauss--Seidel algorithm for solving linear equations and prove its convergence theory. More specifically, we show the exponential convergence of the  norms of the expected iterates via exact formulas (see Theorems \ref{main1} and \ref{main2}) for general (consistent or inconsistent) linear systems, and prove the exponential convergence of the expected norms of the error and the residual (see Theorems \ref{main3} and \ref{main4}) for consistent linear systems.

{\it Organization of this paper}. In the rest of this section, we give some notation. In Section 2 we describe the doubly stochastic block Gauss--Seidel algorithm and prove its convergence theory. In Section 3 we report the numerical results. Finally, we present brief concluding remarks in Section 4.

{\it Notation}. For any random variables $\bm\xi$ and $\bm\zeta$, we use $\mbbe[\bm\xi]$ and $\mbbe[\bm\xi\ |\bm\zeta]$ to denote the expectation of $\bm\xi$ and the conditional expectation of $\bm\xi$ given $\bm\zeta$, respectively. For an integer $m\geq 1$, let $[m]:=\{1,2,3,\ldots,m\}$. Lowercase (upper-case) boldface letters are reserved for column vectors (matrices). For any vector $\mbf u\in\mbbr^m$, we use $\mbf u_i$, $\bf u^\rmt$ and $\|\mbf u\|_2$ to denote, the $i$th entry, the transpose and the Euclidean norm of $\mbf u$, respectively. We use $\mbf I$ to denote the identity matrix whose order is clear from the context. For any matrix $\mbf A\in\mbbr^{m\times n}$, we use $\mbf A_{i,j}$, $\mbf A_{i,:}$, $\mbf A_{:,j}$ $\mbf A^\rmt$, $\mbf A^\dag$, $\|\mbf A\|_\rmf$, $\ran(\mbf A)$, $\rank(\mbf A)$, $\sigma_{1}(\mbf A)\geq\sigma_{2}(\mbf A)\geq\cdots\geq\sigma_{r}(\mbf A)>0$ to denote the $(i,j)$ entry, the $i$th row, the $j$th column, the transpose, the Moore-Penrose pseudoinverse, the Frobenius norm, the column space, the rank, and all the nonzero singular values of $\mbf A$, respectively. Obviously, $\rank(\mbf A)=r$. We call a matrix $\mbf A\in\mbbr^{m\times n}$ full column rank if $\rank(\mbf A)=n$ and rank deficient if $\rank(\mbf A)<n$. For index sets $\mcali\subseteq[m]$ and $\mcalj\subseteq[n]$, let $\mbf A_{\mcali,:}$, $\mbf A_{:,\mcalj}$, and $\mbf A_{\mcali,\mcalj}$ denote the row submatrix indexed by $\mcali$, the column submatrix indexed by $\mcalj$, and the submatrix that lies in the rows indexed by $\mcali$ and the columns indexed by $\mcalj$, respectively. The linear system (\ref{lineq}) is called consistent if $\mbf b\in\ran(\mbf A)$, i.e., a solution exists; otherwise, it is called inconsistent. 

\section{A doubly stochastic block Gauss--Seidel algorithm}

Let $\{\mcali_1,\mcali_2,\ldots,\mcali_s\}$ denote a partition of $[m]$ such that, for $i,j=1,2,\ldots,s$ and $i\neq j,$ $$ \mcali_i\neq\emptyset,\quad\mcali_i\cap\mcali_j=\emptyset,\quad \bigcup_{i=1}^s\mcali_i=[m].$$ Let $\{\mcalj_1,\mcalj_2,\ldots,\mcalj_t\}$ denote a partition of $[n]$ such that, for $i,j=1,2,\ldots,t$ and $i\neq j,$ $$ \mcalj_i\neq\emptyset,\quad\mcalj_i\cap\mcalj_j=\emptyset,\quad \bigcup_{i=1}^t\mcalj_i=[n].$$ Let $$\mcalp=\{\mcali_1,\mcali_2,\ldots,\mcali_s\}\times \{\mcalj_1,\mcalj_2,\ldots,\mcalj_t\}.$$ We propose the following doubly stochastic block Gauss--Seidel algorithm (Algorithm 1) for solving the linear system  $\bf Ax=b$. 

\begin{center}
\begin{tabular*}{132mm}{l}
\toprule {\bf Algorithm 1:} A doubly stochastic block Gauss--Seidel algorithm\\ 
\hline
\qquad Let $\alpha>0$. Initialize $\mbf x^0\in\mbbr^n$\\
\qquad {\bf for} $k=1,2,\ldots,$ {\bf do}\\
\qquad \qquad  Pick $(\mcali,\mcalj)\in\mcalp$ with probability $\dsp\frac{\|\mbf A_{\mcali,\mcalj}\|^2_\rmf}{\|\mbf A\|_\rmf^2}$\\
\qquad \qquad  Set $\dsp\mbf x^k=\mbf x^{k-1}-\alpha \frac{\mbf I_{:,\mcalj}(\mbf A_{\mcali,\mcalj})^\rmt(\mbf I_{:,\mcali})^\rmt}{\|\mbf A_{\mcali,\mcalj}\|_\rmf^2} (\mbf A\mbf x^{k-1}-\mbf b)$\\
\bottomrule
\end{tabular*}
\end{center}

Here we consider constant step size for simplicity.  By varying the row partition parameter $s$ and the column partition parameter $t$, we recover the following well-known algorithms as special cases: \bit 
\item[$\bullet$] Landweber \cite{landweber1951itera} ($s=1$ and $t=1$), $$\mbf x^k=\mbf x^{k-1}-\alpha\frac{{\bf A}^\rmt}{\|\mbf A\|_\rmf^2}(\mbf A\mbf x^{k-1}-\mbf b).$$
\item[$\bullet$] Randomized Kaczmarz \cite{strohmer2009rando} ($s=m$ and $t=1$), $$\mbf x^k=\mbf x^{k-1}-\alpha\frac{\mbf A_{i,:}\mbf x^{k-1}-\mbf b_i}{\|\mbf A_{i,:}\|_2^2}(\mbf A_{i,:})^\rmt.$$
\item[$\bullet$] Randomized Gauss--Seidel \cite{leventhal2010rando} ($s=1$ and $t=n$), $$\dsp\mbf x^k=\mbf x^{k-1}-\alpha\frac{(\mbf A_{:,j})^\rmt(\mbf A\mbf x^{k-1}-\mbf b)}{\|\mbf A_{:,j}\|_2^2} \mbf I_{:,j}.$$ 
\item[$\bullet$] Doubly Stochastic Gauss--Seidel \cite{razaviyayn2019linea} ($s=m$ and $t=n$), $$\dsp\mbf x^k=\mbf x^{k-1}-\alpha \frac{\mbf A_{i,j}(\mbf A_{i,:}\mbf x^{k-1}-\mbf b_i)}{|\mbf A_{i,j}|^2}\mbf I_{:,j}.$$
\eit

The conditional expectation of $\mbf x^k$ given $\mbf x^{k-1}$ is \beqas\mbbe[\mbf x^k\ |\mbf x^{k-1}]&=&\mbf x^{k-1}-\alpha\mbbe\bem\dsp\frac{\mbf I_{:,\mcalj}(\mbf A_{\mcali,\mcalj})^\rmt(\mbf I_{:,\mcali})^\rmt}{\|\mbf A_{\mcali,\mcalj}\|_\rmf^2}\eem(\mbf A\mbf x^{k-1}-\mbf b)\\
&=&\mbf x^{k-1}-\alpha\l(\sum_{(\mcali,\mcalj)\in\mcalp}\frac{\mbf I_{:,\mcalj}(\mbf A_{\mcali,\mcalj})^\rmt(\mbf I_{:,\mcali})^\rmt}{\|\mbf A_{\mcali,\mcalj}\|_\rmf^2}\frac{\|\mbf A_{\mcali,\mcalj}\|_\rmf^2}{\|\mbf A\|_\rmf^2}\r)(\mbf A\mbf x^{k-1}-\mbf b)\\
&=&\mbf x^{k-1}-\alpha\frac{{\bf A}^\rmt}{\|\mbf A\|_\rmf^2}(\mbf A\mbf x^{k-1}-\mbf b). \eeqas
Note that the gradient of the objective function of the optimization problem (\ref{opt}) is $$\nabla f(\mbf x)=\frac{\mbf A^\rmt}{\|\mbf A\|_\rmf^2}(\mbf A\mbf x-\mbf b).$$ It follows 
$$\mbbe[\mbf x^k\ |\mbf x^{k-1}]=\mbf x^{k-1}-\alpha\nabla f(\mbf x^{k-1}).$$
Therefore, DSBGS can be viewed as a stochastic gradient descent method for solving the optimization problem (\ref{opt}).

\subsection{The exponential convergence of the norms of the expected iterates}
In this subsection we show the exponential convergence of the norms of the expected iterates for general (consistent or inconsistent) linear systems. The following lemma will be used to prove Theorems \ref{main1} and \ref{main2}. Its proof (via singular value decomposition) is straightforward and we omit the details.

\begin{lemma}\label{leq} Let $\alpha>0$ and $\mbf A$ be any nonzero real matrix. For every $\mbf u\in\ran(\mbf A)$, it holds $$\l\|\l(\mbf I-\frac{\alpha\bf AA^\rmt}{\|\mbf A\|_\rmf^2}\r)^k\mbf u\r\|_2\leq\l(\max_{1\leq i\leq r}\l|1-\frac{\alpha\sigma_i^2(\mbf A)}{\|\mbf A\|_\rmf^2}\r|\r)^k\|\mbf u\|_2.$$ 
\end{lemma}


In Theorem \ref{main1}, we show the exponential convergence of the norm of the expected error for {\it consistent} linear systems.
 
\begin{theorem}\label{main1} Let $\mbf x^k$ denote the $k$th iterate of {\rm DSBGS} applied to the consistent linear system  $\bf Ax=b$ with arbitrary $\mbf x^0\in\mbbr^n$. 
It holds $$\|\mbbe[{\bf x}^k-\mbf x^0_\star]\|_2\leq\l(\max_{1\leq i\leq r}\l|1-\frac{\alpha\sigma_i^2(\mbf A)}{\|\mbf A\|_\rmf^2}\r|\r)^k\|{\bf x}^0-\mbf x^0_\star\|_2,$$ where $$\mbf x^0_\star:=(\mbf I-\mbf A^\dag\mbf A)\mbf x^0+\mbf A^\dag\mbf b$$ is the orthogonal projection of $\mbf x^0$ onto the solution set $\{\mbf x\in\mbbr^n\ |\ \bf  A x= b\}$.  
\end{theorem}

\proof The conditional expectation of $\mbf x^k-\mbf x_\star^0$ given $\mbf x^{k-1}$ is \beqas \mbbe[{\bf x}^k-\mbf x_\star^0\ |\mbf x^{k-1}]&=&\mbbe[\mbf x^k \ |\mbf x^{k-1}] -\mbf x_\star^0\\ 
&=&\mbf x^{k-1}-\alpha\frac{{\bf A}^\rmt}{\|\mbf A\|_\rmf^2}(\mbf A\mbf x^{k-1}-\mbf b)-\mbf x_\star^0\\&=& \mbf x^{k-1}-\alpha\frac{{\bf A}^\rmt}{\|\mbf A\|_\rmf^2}(\mbf A\mbf x^{k-1}-\mbf A\mbf x_\star^0)-\mbf x_\star^0\\
&=&\l(\mbf I-\frac{\alpha\mbf A^\rmt\mbf A}{\|\mbf A\|_\rmf^2}\r) ({\bf x}^{k-1}-\mbf x_\star^0).\eeqas
Taking expectation gives \beqas\mbbe[{\bf x}^k-\mbf x_\star^0]&=&\mbbe[\mbbe[{\bf x}^k-\mbf x_\star^0\ |\mbf x^{k-1}]]\\ &=& \l(\mbf I-\frac{\alpha\mbf A^\rmt\mbf A}{\|\mbf A\|_\rmf^2}\r) \mbbe[{\bf x}^{k-1}-\mbf x_\star^0]\\ &=& \l(\mbf I-\frac{\alpha\mbf A^\rmt\mbf A}{\|\mbf A\|_\rmf^2}\r)^k ({\bf x}^0-\mbf x_\star^0). \eeqas 
Applying the norms to both sides we obtain $$\|\mbbe[{\bf x}^k-\mbf x_\star^0]\|_2\leq\l(\max_{1\leq i\leq r}\l|1-\frac{\alpha\sigma_i^2(\mbf A)}{\|\mbf A\|_\rmf^2}\r|\r)^k\|{\bf x}^0-\mbf x_\star^0\|_2.$$ Here the inequality follows from the fact that $${\bf x}^0-\mbf x_\star^0={\bf A^\dag Ax}^0-\mbf A^\dag\mbf b\in\ran(\mbf A^\rmt)$$ and Lemma \ref{leq}.
\qed

\begin{remark} If $\mbf x^0\in\ran(\mbf A^\rmt)$, then $\mbf x^0_\star=\mbf A^\dag\mbf b$.
\end{remark}

\begin{remark} To ensure convergence of the expected iterate, it suffices to have $$\dsp\max_{1\leq i\leq r}\l|1-\frac{\alpha\sigma_i^2(\mbf A)}{\|\mbf A\|_\rmf^2}\r|<1\quad i.e.,\quad  0<\alpha<\dsp\frac{2\|\mbf A\|_\rmf^2}{\sigma_1^2(\mbf A)}.$$ 
\end{remark}

In Theorem \ref{main2}, we show the exponential convergence of  $\|\mbbe[{\bf Ax}^k-{\bf Ax}_\star]\|_2$ for the {\it consistent or inconsistent} linear system $\bf Ax=b$, where $\mbf x_\star$ is any solution of the normal equations $$\mbf A^\rmt\mbf A\mbf x=\mbf A^\rmt\mbf b.$$

\begin{theorem}\label{main2} Let $\mbf x^k$ denote the $k$th iterate of {\rm DSBGS} applied to the consistent or inconsistent linear system  $\bf Ax=b$ with arbitrary $\mbf x^0\in\mbbr^n$. It holds $$\|\mbbe[{\bf Ax}^k-{\bf Ax}_\star]\|_2\leq\l(\max_{1\leq i\leq r}\l|1-\frac{\alpha\sigma_i^2(\mbf A)}{\|\mbf A\|_\rmf^2}\r|\r)^k\|{\bf Ax}^0-{\bf Ax}_\star\|_2,$$ where $\mbf x_\star$ is any solution of $\mbf A^\rmt\mbf A\mbf x=\mbf A^\rmt\mbf b$.
\end{theorem}

\proof The conditional expectation of ${\bf Ax}^k-{\bf Ax}_\star$ given $\mbf x^{k-1}$ is \beqas \mbbe[{\bf Ax}^k-{\bf Ax}_\star\ |\mbf x^{k-1}]&=&\mbf A(\mbbe[\mbf x^k \ |\mbf x^{k-1}] - {\bf x_\star})\\ &=& \mbf A\l(\mbf x^{k-1}-\alpha\frac{{\bf A}^\rmt}{\|\mbf A\|_\rmf^2}(\mbf A\mbf x^{k-1}-\mbf b)- {\bf x_\star}\r)  \\ &=& \mbf A\l(\mbf x^{k-1}-\frac{\alpha\mbf A^\rmt}{\|\mbf A\|_\rmf^2}(\mbf A\mbf x^{k-1}-\mbf A\mbf x_\star)- {\bf x_\star}\r) \qquad{\rm (by\  \bf A^\rmt b=\bf A^\rmt  Ax_\star)} \\ &=& \mbf A\mbf x^{k-1} - {\bf Ax_\star}-\frac{\alpha\mbf A\mbf A^\rmt}{\|\mbf A\|_\rmf^2}(\mbf A\mbf x^{k-1}-\mbf A\mbf x_\star)\\ &=&\l(\mbf I-\frac{\alpha\mbf A\mbf A^\rmt}{\|\mbf A\|_\rmf^2}\r) ({\bf Ax}^{k-1}-{\bf Ax_\star}).\eeqas
Taking expectation gives \beqas\mbbe[{\bf Ax}^k-{\bf Ax}_\star]&=&\mbbe[\mbbe[{\bf Ax}^k-{\bf Ax}_\star\ |\mbf x^{k-1}]]\\ &=& \l(\mbf I-\frac{\alpha\mbf A\mbf A^\rmt}{\|\mbf A\|_\rmf^2}\r) \mbbe[{\bf Ax}^{k-1}-{\bf Ax}_\star]\\ &=& \l(\mbf I-\frac{\alpha\mbf A\mbf A^\rmt}{\|\mbf A\|_\rmf^2}\r)^k ({\bf Ax}^0-{\bf Ax}_\star). \eeqas 
Applying the norms to both sides we obtain $$\|\mbbe[{\bf Ax}^k-{\bf Ax}_\star]\|_2\leq\l(\max_{1\leq i\leq r}\l|1-\frac{\alpha\sigma_i^2(\mbf A)}{\|\mbf A\|_\rmf^2}\r|\r)^k\|{\bf Ax}^0-{\bf Ax}_\star\|_2.$$ Here the inequality follows from the fact that ${\bf Ax}^0-{\bf Ax}_\star\in\ran(\mbf A)$ and Lemma \ref{leq}.
\qed

\subsection{The exponential convergence of the expected norms of the error and the residual}

In this subsection we prove the exponential convergence of the expected norms of the error or the residual for consistent linear systems. {The convergence depends on the positive number $\rho$  defined as $$\beta=\max_{(\mcali,\mcalj)\in\mcalp}\frac{\|\mbf A_{\mcali,\mcalj}\|_2^2}{\|\mbf A_{\mcali,\mcalj}\|_\rmf^2}.$$} The following two lemmas will be used. Their proofs are straightforward and we omit the details.
\begin{lemma}\label{lemma2} For any vector $\mbf u\in\mbbr^m$ and any matrix $\mbf A\in\mbbr^{m\times n}$, it holds $${\mbf u^\rmt\mbf A\mbf A^\rmt\mbf u\leq\|\mbf A\|_2^2\mbf u^\rmt\mbf u.}$$
\end{lemma}
\begin{lemma}\label{lemma3} For any matrix $\mbf A\in\mbbr^{m\times n}$ with rank $r$ and any vector $\mbf u\in\ran(\mbf A)$, it holds $$\mbf u^\rmt\mbf A\mbf A^\rmt\mbf u\geq\sigma_r^2({\mbf A})\|\mbf u\|_2^2.$$
\end{lemma}

For {\it full column rank  consistent} linear systems, we prove the exponential convergence of the expected norm of the error in the following theorem. We recall that in this case  $\bf A^\dag b$ is the unique solution of $\bf Ax=b$.

\begin{theorem}\label{main3} Let $\mbf x^k$ denote the $k$th iterate of {\rm DSBGS} applied to the full column rank consistent linear system $\bf Ax=b$ with arbitrary $\mbf x^0\in\mbbr^n$. Assume  $0<\alpha<2/{(t\beta)}$. It holds $$\mbbe[\|{\bf x}^k-\mbf A^\dag\mbf b\|_2^2]\leq\l(1-\frac{(2\alpha-t{\beta}\alpha^2)\sigma_n^2(\mbf A)}{\|\mbf A\|_\rmf^2}\r)^k \|\mbf x^0-\mbf A^\dag\mbf b\|_2^2.$$
\end{theorem}
\proof Note that
\beqas \|\mbf x^k-\mbf A^\dag\mbf b\|_2^2&=&\l\|\mbf x^{k-1}-\alpha\l(\frac{\mbf I_{:,\mcalj}(\mbf A_{\mcali,\mcalj})^\rmt(\mbf I_{:,\mcali})^\rmt}{\|\mbf A_{\mcali,\mcalj}\|_\rmf^2}\r)(\mbf A\mbf x^{k-1}-\mbf b)-\mbf A^\dag\mbf b\r\|_2^2\\
&=& \l\|\mbf x^{k-1}-\alpha\l(\frac{\mbf I_{:,\mcalj}(\mbf A_{\mcali,\mcalj})^\rmt(\mbf I_{:,\mcali})^\rmt}{\|\mbf A_{\mcali,\mcalj}\|_\rmf^2}\r)\mbf A(\mbf x^{k-1}-\mbf A^\dag\mbf b)-\mbf A^\dag\mbf b\r\|_2^2\\
&=& \l\|\mbf x^{k-1}-\mbf A^\dag\mbf b-\alpha\l(\frac{\mbf I_{:,\mcalj}(\mbf A_{\mcali,\mcalj})^\rmt(\mbf I_{:,\mcali})^\rmt\mbf A}{\|\mbf A_{\mcali,\mcalj}\|_\rmf^2}\r)(\mbf x^{k-1}-\mbf A^\dag\mbf b)\r\|_2^2 \\
&=& \|\mbf x^{k-1}-\mbf A^\dag\mbf b\|_2^2-2\alpha(\mbf x^{k-1}-\mbf A^\dag\mbf b)^\rmt\l(\frac{\mbf I_{:,\mcalj}(\mbf A_{\mcali,\mcalj})^\rmt(\mbf I_{:,\mcali})^\rmt\mbf A}{\|\mbf A_{\mcali,\mcalj}\|_\rmf^2}\r)(\mbf x^{k-1}-\mbf A^\dag\mbf b)\\
&& + \alpha^2(\mbf x^{k-1}-\mbf A^\dag\mbf b)^\rmt\l(\frac{\mbf A^\rmt\mbf I_{:,\mcali}\mbf A_{\mcali,\mcalj}(\mbf A_{\mcali,\mcalj})^\rmt(\mbf I_{:,\mcali})^\rmt\mbf A}{\|\mbf A_{\mcali,\mcalj}\|_\rmf^4}\r)(\mbf x^{k-1}-\mbf A^\dag\mbf b)\\ 
&\leq&  \|\mbf x^{k-1}-\mbf A^\dag\mbf b\|_2^2-2\alpha(\mbf x^{k-1}-\mbf A^\dag\mbf b)^\rmt\l(\frac{\mbf I_{:,\mcalj}(\mbf A_{\mcali,\mcalj})^\rmt(\mbf I_{:,\mcali})^\rmt\mbf A}{\|\mbf A_{\mcali,\mcalj}\|_\rmf^2}\r)(\mbf x^{k-1}-\mbf A^\dag\mbf b)\\
&& + {\alpha^2(\mbf x^{k-1}-\mbf A^\dag\mbf b)^\rmt\l(\frac{\|\mbf A_{\mcali,\mcalj}\|_2^2}{\|\mbf A_{\mcali,\mcalj}\|_\rmf^2}\cdot\frac{\mbf A^\rmt\mbf I_{:,\mcali}(\mbf I_{:,\mcali})^\rmt\mbf A}{\|\mbf A_{\mcali,\mcalj}\|_\rmf^2}\r)(\mbf x^{k-1}-\mbf A^\dag\mbf b)\quad \mbox{(by Lemma \ref{lemma2})} }\\
& \leq & { \|\mbf x^{k-1}-\mbf A^\dag\mbf b\|_2^2-2\alpha(\mbf x^{k-1}-\mbf A^\dag\mbf b)^\rmt\l(\frac{\mbf I_{:,\mcalj}(\mbf A_{\mcali,\mcalj})^\rmt(\mbf I_{:,\mcali})^\rmt\mbf A}{\|\mbf A_{\mcali,\mcalj}\|_\rmf^2}\r)(\mbf x^{k-1}-\mbf A^\dag\mbf b)}\\
&&  {+ \beta\alpha^2(\mbf x^{k-1}-\mbf A^\dag\mbf b)^\rmt\l(\frac{\mbf A^\rmt\mbf I_{:,\mcali}(\mbf I_{:,\mcali})^\rmt\mbf A}{\|\mbf A_{\mcali,\mcalj}\|_\rmf^2}\r)(\mbf x^{k-1}-\mbf A^\dag\mbf b)}.
\eeqas
Taking expectation gives 
\beqas
\mbbe[\|{\bf x}^k-\mbf A^\dag\mbf b\|_2^2\ |\mbf x^{k-1}]&\leq&\|\mbf x^{k-1}-\mbf A^\dag\mbf b\|_2^2-(2\alpha-t{\beta}\alpha^2)(\mbf x^{k-1}-\mbf A^\dag\mbf b)^\rmt\l(\frac{\mbf A ^\rmt\mbf A}{\|\mbf A\|_\rmf^2}\r)(\mbf x^{k-1}-\mbf A^\dag\mbf b)\\
&\leq&\l(1-\frac{(2\alpha-t{\beta}\alpha^2)\sigma_n^2(\mbf A)}{\|\mbf A\|_\rmf^2}\r)\|\mbf x^{k-1}-\mbf A^\dag\mbf b\|_2^2.\quad ({\rm by\ Lemma\ \ref{lemma3}})
\eeqas
Taking expectation again gives  
\beqas
\mbbe[\|{\bf x}^k-\mbf A^\dag\mbf b\|_2^2]&=&\mbbe[\mbbe[\|{\bf x}^k-\mbf A^\dag\mbf b\|_2^2\ |\mbf x^{k-1}]]
\\ &\leq&\l(1-\frac{(2\alpha-t{\beta}\alpha^2)\sigma_n^2(\mbf A)}{\|\mbf A\|_\rmf^2}\r)\mbbe[\|\mbf x^{k-1}-\mbf A^\dag\mbf b\|_2^2]\\
&\leq&\l(1-\frac{(2\alpha-t{\beta}\alpha^2)\sigma_n^2(\mbf A)}{\|\mbf A\|_\rmf^2}\r)^k \|\mbf x^0-\mbf A^\dag\mbf b\|_2^2.\qquad \qed\eeqas

\begin{remark}
For the special case $s=m$, $t=1$, $\alpha=1$  (i.e., the randomized Kaczmarz algorithm{, we have $\beta=1$ in this case}) and the special case $s=m$, $t=n$, $\alpha=1/n$ (i.e., the doubly stochastic Gauss--Seidel algorithm{,  we have $\beta=1$ in this case}), the results of Theorem \ref{main3} are given in \cite[Theorem 2]{strohmer2009rando} and \cite[Theorem 1]{razaviyayn2019linea}, respectively.
\end{remark}

\begin{remark} For rank deficient  consistent linear systems, if $t=1$ and $\mbf x^0\in\mbbr^n$, we can show $\mbf x^k-\mbf x_\star^0\in\ran(\mbf A^\rmt)$ by induction, where $\mbf x^0_\star=(\mbf I-\mbf A^\dag\mbf A)\mbf x^0+\mbf A^\dag\mbf b$ is the orthogonal projection of $\mbf x^0$ onto the solution set $\{\mbf x\in\mbbr^n\ |\ \bf  A x= b\}$.  Then by the same approach as used in the proof of Theorem \ref{main3}, for any $s\in[m]$ and $t=1$, we can prove the convergence bound  $$\mbbe[\|{\bf x}^k-\mbf x_\star^0\|_2^2]\leq\l(1-\frac{(2\alpha-{\beta}\alpha^2)\sigma_r^2(\mbf A)}{\|\mbf A\|_\rmf^2}\r)^k \|\mbf x^0-\mbf x_\star^0\|_2^2.$$ This result for the special case $s=m$ and $t=1$ (i.e., the randomized Kaczmarz algorithm{, we have $\beta=1$ in this case}) with $\mbf x^0\in\ran(\mbf A^\rmt) $ is  given in \cite[Theorem 3.4]{zouzias2013rando}.
\end{remark}

Next, we prove the exponential convergence of the expected norm of the residual for {\it full column rank or rank-deficient consistent} linear systems.

\begin{theorem}\label{main4} Let $\mbf x^k$ denote the $k$th iterate of {\rm DSBGS} applied to the consistent linear system (full column rank or rank deficient) $\bf Ax=b$ with arbitrary $\mbf x^0\in\mbbr^n$. If $t=n$ and $0<\alpha<2\sigma_r^2(\mbf A)/{(\beta\|\mbf A\|_\rmf^2)}$, then  $$\mbbe[\|\mbf A\mbf x^k-\mbf b\|_2^2]\leq\l(1+{\beta}\alpha^2-\frac{2\alpha\sigma_r^2(\mbf A)}{\|\mbf A\|_\rmf^2}\r)^k\|\mbf A\mbf x^0-\mbf b\|_2^2.$$ If $t<n$ and $0<\alpha<2\sigma_r^2(\mbf A)/(t\rho{\beta})$, then $$\mbbe[\|\mbf A\mbf x^k-\mbf b\|_2^2]\leq\l(1-\frac{2\alpha\sigma_r^2(\mbf A)-t\rho{\beta}\alpha^2}{\|\mbf A\|_\rmf^2}\r)^k\|\mbf A\mbf x^0-\mbf b\|_2^2,$$ where $$\rho=\max_{1\leq j\leq t}\sigma_1^2(\mbf A_{:,\mcalj_j}).$$
\end{theorem}

\proof Note that
\beqas \|{\bf Ax}^k-\mbf b\|_2^2&=&\l\|{\bf Ax}^{k-1}-\alpha\l(\frac{\mbf A\mbf I_{:,\mcalj}(\mbf A_{\mcali,\mcalj})^\rmt(\mbf I_{:,\mcali})^\rmt}{\|\mbf A_{\mcali,\mcalj}\|_\rmf^2}\r)(\mbf A\mbf x^{k-1}-\mbf b)-\mbf b\r\|_2^2\\
&=& \|{\bf Ax}^{k-1}-\mbf b\|_2^2-2\alpha({\bf Ax}^{k-1}-\mbf b)^\rmt\l(\frac{\mbf A\mbf I_{:,\mcalj}(\mbf A_{\mcali,\mcalj})^\rmt(\mbf I_{:,\mcali})^\rmt}{\|\mbf A_{\mcali,\mcalj}\|_\rmf^2}\r)({\bf Ax}^{k-1}-\mbf b)\\
&& + \alpha^2({\bf Ax}^{k-1}-\mbf b)^\rmt\l(\frac{\mbf I_{:,\mcali}\mbf A_{\mcali,\mcalj}(\mbf I_{:,\mcalj})^\rmt\mbf A^\rmt\mbf A \mbf I_{:,\mcalj}(\mbf A_{\mcali,\mcalj})^\rmt(\mbf I_{:,\mcali})^\rmt}{\|\mbf A_{\mcali,\mcalj}\|_\rmf^4}\r)({\bf Ax}^{k-1}-\mbf b). \eeqas
If $t=n$, then it follows from $(\mbf I_{:,\mcalj})^\rmt\mbf A^\rmt\mbf A \mbf I_{:,\mcalj}=\|\mbf A_{:,\mcalj}\|_\rmf^2$ (since $\mbf A \mbf I_{:,\mcalj}=\mbf A_{:,\mcalj}$ is a column vector) that
\beqas \|{\bf Ax}^k-\mbf b\|_2^2 &=&  
\|{\bf Ax}^{k-1}-\mbf b\|_2^2-2\alpha({\bf Ax}^{k-1}-\mbf b)^\rmt\l(\frac{\mbf A\mbf I_{:,\mcalj}(\mbf A_{\mcali,\mcalj})^\rmt(\mbf I_{:,\mcali})^\rmt}{\|\mbf A_{\mcali,\mcalj}\|_\rmf^2}\r)({\bf Ax}^{k-1}-\mbf b)\\
&& + \alpha^2({\bf Ax}^{k-1}-\mbf b)^\rmt\l(\frac{\|\mbf A_{:,\mcalj}\|_\rmf^2\mbf I_{:,\mcali}\mbf A_{\mcali,\mcalj}(\mbf A_{\mcali,\mcalj})^\rmt(\mbf I_{:,\mcali})^\rmt}{\|\mbf A_{\mcali,\mcalj}\|_\rmf^4}\r)({\bf Ax}^{k-1}-\mbf b)\\
&\leq &
\|{\bf Ax}^{k-1}-\mbf b\|_2^2-2\alpha({\bf Ax}^{k-1}-\mbf b)^\rmt\l(\frac{\mbf A\mbf I_{:,\mcalj}(\mbf A_{\mcali,\mcalj})^\rmt(\mbf I_{:,\mcali})^\rmt}{\|\mbf A_{\mcali,\mcalj}\|_\rmf^2}\r)({\bf Ax}^{k-1}-\mbf b)\\
&& + \alpha^2({\bf Ax}^{k-1}-\mbf b)^\rmt\l({\frac{\|\mbf A_{\mcali,\mcalj}\|_2^2}{\|\mbf A_{\mcali,\mcalj}\|_\rmf^2}}\frac{\|\mbf A_{:,\mcalj}\|_\rmf^2\mbf I_{:,\mcali}(\mbf I_{:,\mcali})^\rmt}{\|\mbf A_{\mcali,\mcalj}\|_\rmf^2}\r)({\bf Ax}^{k-1}-\mbf b) \quad ({\rm by\ Lemma\ \ref{lemma2}})\\
&\leq &
{\|{\bf Ax}^{k-1}-\mbf b\|_2^2-2\alpha({\bf Ax}^{k-1}-\mbf b)^\rmt\l(\frac{\mbf A\mbf I_{:,\mcalj}(\mbf A_{\mcali,\mcalj})^\rmt(\mbf I_{:,\mcali})^\rmt}{\|\mbf A_{\mcali,\mcalj}\|_\rmf^2}\r)({\bf Ax}^{k-1}-\mbf b)}\\
&& {+ \beta\alpha^2({\bf Ax}^{k-1}-\mbf b)^\rmt\l(\frac{\|\mbf A_{:,\mcalj}\|_\rmf^2\mbf I_{:,\mcali}(\mbf I_{:,\mcali})^\rmt}{\|\mbf A_{\mcali,\mcalj}\|_\rmf^2}\r)({\bf Ax}^{k-1}-\mbf b)}
\eeqas
Taking expectation gives
\beqas\mbbe[\|{\bf Ax}^k-\mbf b\|_2^2\ |\mbf x^{k-1}] &\leq&(1+{\beta}\alpha^2)\|{\bf Ax}^{k-1}-\mbf b\|_2^2-2\alpha({\bf Ax}^{k-1}-\mbf b)^\rmt\l(\frac{\mbf A\mbf A^\rmt}{\|\mbf A \|_\rmf^2}\r)({\bf Ax}^{k-1}-\mbf b)\\ &\leq&
\l(1+{\beta}\alpha^2-\frac{2\alpha\sigma_r^2(\mbf A)}{\|\mbf A\|_\rmf^2}\r)\|\mbf A\mbf x^{k-1}-\mbf b\|_2^2.\eeqas 
The last inequality follows from $\mbf A\mbf x^{k-1}-\mbf b\in\ran(\mbf A)$ and Lemma \ref{lemma3}.
Taking expectation again gives \beqas \mbbe[\|{\bf Ax}^k-\mbf b\|_2^2]&=&\mbbe[ \mbbe[\|{\bf Ax}^k-\mbf b\|_2^2\ |\mbf x^{k-1}]]\\
&\leq& \l(1+{\beta}\alpha^2-\frac{2\alpha\sigma_r^2(\mbf A)}{\|\mbf A\|_\rmf^2}\r)\mbbe[\|\mbf A\mbf x^{k-1}-\mbf b\|_2^2]\\
&\leq& \l(1+{\beta}\alpha^2-\frac{2\alpha\sigma_r^2(\mbf A)}{\|\mbf A\|_\rmf^2}\r)^k\|\mbf A\mbf x^0-\mbf b\|_2^2.\eeqas
If $t<n$, then it follows from $(\mbf I_{:,\mcalj})^\rmt\mbf A^\rmt\mbf A \mbf I_{:,\mcalj}=\mbf A_{:,\mcalj}^\rmt\mbf A_{:,\mcalj}\preceq\rho\mbf I$ (since $\rho=\max_{1\leq j\leq t}\sigma_1^2(\mbf A_{:,\mcalj_j})$) that
\beqas \|{\bf Ax}^k-\mbf b\|_2^2 &\leq&  
\|{\bf Ax}^{k-1}-\mbf b\|_2^2-2\alpha({\bf Ax}^{k-1}-\mbf b)^\rmt\l(\frac{\mbf A\mbf I_{:,\mcalj}(\mbf A_{\mcali,\mcalj})^\rmt(\mbf I_{:,\mcali})^\rmt}{\|\mbf A_{\mcali,\mcalj}\|_\rmf^2}\r)({\bf Ax}^{k-1}-\mbf b)\\
&& + \alpha^2({\bf Ax}^{k-1}-\mbf b)^\rmt\l(\frac{\rho\mbf I_{:,\mcali}\mbf A_{\mcali,\mcalj}(\mbf A_{\mcali,\mcalj})^\rmt(\mbf I_{:,\mcali})^\rmt}{\|\mbf A_{\mcali,\mcalj}\|_\rmf^4}\r)({\bf Ax}^{k-1}-\mbf b)\\
&\leq&
\|{\bf Ax}^{k-1}-\mbf b\|_2^2-2\alpha({\bf Ax}^{k-1}-\mbf b)^\rmt\l(\frac{\mbf A\mbf I_{:,\mcalj}(\mbf A_{\mcali,\mcalj})^\rmt(\mbf I_{:,\mcali})^\rmt}{\|\mbf A_{\mcali,\mcalj}\|_\rmf^2}\r)({\bf Ax}^{k-1}-\mbf b)\\
&& + \alpha^2({\bf Ax}^{k-1}-\mbf b)^\rmt\l({\frac{\|\mbf A_{\mcali,\mcalj}\|_2^2}{\|\mbf A_{\mcali,\mcalj}\|_\rmf^2}}\frac{\rho\mbf I_{:,\mcali}(\mbf I_{:,\mcali})^\rmt}{\|\mbf A_{\mcali,\mcalj}\|_\rmf^2}\r)({\bf Ax}^{k-1}-\mbf b) \quad ({\rm by\ Lemma\ \ref{lemma2}})\\
&\leq&
\|{\bf Ax}^{k-1}-\mbf b\|_2^2-2\alpha({\bf Ax}^{k-1}-\mbf b)^\rmt\l(\frac{\mbf A\mbf I_{:,\mcalj}(\mbf A_{\mcali,\mcalj})^\rmt(\mbf I_{:,\mcali})^\rmt}{\|\mbf A_{\mcali,\mcalj}\|_\rmf^2}\r)({\bf Ax}^{k-1}-\mbf b)\\
&& + \alpha^2({\bf Ax}^{k-1}-\mbf b)^\rmt\l(\frac{\rho{\beta}\mbf I_{:,\mcali}(\mbf I_{:,\mcali})^\rmt}{\|\mbf A_{\mcali,\mcalj}\|_\rmf^2}\r)({\bf Ax}^{k-1}-\mbf b).
\eeqas
Taking expectation gives
\beqas\mbbe[\|{\bf Ax}^k-\mbf b\|_2^2\ |\mbf x^{k-1}] &\leq&\l(1+\frac{t\rho{\beta}\alpha^2}{\|\mbf A\|_\rmf^2}\r)\|{\bf Ax}^{k-1}-\mbf b\|_2^2-2\alpha({\bf Ax}^{k-1}-\mbf b)^\rmt\l(\frac{\mbf A\mbf A^\rmt}{\|\mbf A \|_\rmf^2}\r)({\bf Ax}^{k-1}-\mbf b)\\&\leq&
\l(1-\frac{2\alpha\sigma_r^2(\mbf A)-t\rho{\beta}\alpha^2}{\|\mbf A\|_\rmf^2}\r)\|\mbf A\mbf x^{k-1}-\mbf b\|_2^2.\eeqas
The last inequality follows from $\mbf A\mbf x^{k-1}-\mbf b\in\ran(\mbf A)$ and Lemma \ref{lemma3}. Taking expectation again gives \beqas \mbbe[\|{\bf Ax}^k-\mbf b\|_2^2]&=&\mbbe[ \mbbe[\|{\bf Ax}^k-\mbf b\|_2^2\ |\mbf x^{k-1}]]\\
&\leq& \l(1-\frac{2\alpha\sigma_r^2(\mbf A)-t\rho{\beta}\alpha^2}{\|\mbf A\|_\rmf^2}\r)\mbbe[\|\mbf A\mbf x^{k-1}-\mbf b\|_2^2]\\
&\leq& \l(1-\frac{2\alpha\sigma_r^2(\mbf A)-t\rho{\beta}\alpha^2}{\|\mbf A\|_\rmf^2}\r)^k\|\mbf A\mbf x^0-\mbf b\|_2^2.\qquad \qed\eeqas

\begin{remark}
For the special case $s=1$, $t=n$, $\alpha=\sigma_r^2(\mbf A)/\|\mbf A\|_\rmf^2$  (i.e., the randomized Gauss--Seidel algorithm{, we have $\beta=1$ in this case}) and the special case $s=m$, $t=n$, $\alpha=\sigma_r^2(\mbf A)/\|\mbf A\|_\rmf^2$ (i.e., the doubly stochastic Gauss--Seidel algorithm{, we have $\beta=1$ in this case}), the results of Theorem \ref{main4} are given in \cite[Theorem 3.2]{leventhal2010rando} and \cite[Theorem 2]{razaviyayn2019linea}, respectively.
\end{remark}

\begin{remark}  Let $\mbf x_\star$ be any solution of $\bf A^\rmt A x=A^\rmt b$. We have $(\mbf A_{:,\mcalj})^\rmt\mbf b=(\mbf A_{:,\mcalj})^\rmt\mbf A\mbf x_\star$. Then for $s=1$ and $t\in[n]$, by the same approach as used in the proof of Theorem \ref{main4}, we can prove that the convergence bounds (replacing $\mbf b$ by $\mbf A\mbf x_\star$) in Theorem \ref{main4} still hold for {\it inconsistent} linear systems. The result for the  special case $s=1$ and $t=n$ (i.e., the randomized Gauss--Seidel algorithm) was already given in the literature, for example, \cite[Theorem 3.2]{leventhal2010rando}, \cite[Lemma 4.2]{ma2015conve} and \cite[Theorem 3]{du2019tight}. 
\end{remark}

\section{Numerical results}

In this section, we compare the performance of the doubly stochastic block Gauss--Seidel (DSBGS) algorithm proposed in this paper against the randomized Kaczmarz (RK) algorithm for solving consistent linear systems. We only run on small or medium-scale problems. The purpose is to demonstrate that even in these simple examples, DSBGS offers significant advantages over RK. All experiments are performed using MATLAB (R2019a) on a laptop with 2.7-GHz Intel Core i7 processor, 16 GB memory, and macOS Sierra (version 10.12.6).

We use DSBGS($\alpha,\ell,\tau$) to denote the doubly stochastic block Gauss--Seidel algorithm employing the step size $\alpha$, the row partition $\{\mcali_i\}_{i=1}^s$ with $s=\lceil\frac{m}{\ell}\rceil$: \beqas\mcali_i&=&\{(i-1)\ell+1,(i-1)\ell+2,\ldots,i\ell\},\quad i=1,2,\ldots,s-1,\\ \mcali_s&=&\{(s-1)\ell+1,(s-1)\ell+2,\ldots,m\},\eeqas and the column partition $\{\mcalj_j\}_{j=1}^t$ with $t=\lceil\frac{n}{\tau}\rceil$: \beqas\mcalj_j&=&\{(j-1)\tau+1,(j-1)\tau+2,\ldots,j\tau\},\quad j=1,2,\ldots,t-1,\\ \mcalj_t&=&\{(t-1)\tau+1,(t-1)\tau+2,\ldots,n\}.\eeqas The randomized Kaczmarz algorithm with step size $\alpha=1$ is the special case DSBGS($1,1,n$).

To construct a consistent linear system, for a given coefficient matrix $\mbf A$, we set $\bf b = Ax$ where $\mbf x$ is a vector with entries generated from a standard normal distribution. All algorithms are started from the initial guess $\mbf x^0 = \mbf 0$, terminated if $\|\mbf x^k-\mbf A^\dag\mbf b\|_2\leq10^{-5}$.  We report the average number of iterations (denoted as ITER) of RK and DSBGS. We also report the average computing time in seconds (denoted as CPU) and the speed-up of DSBGS against RK, which is defined as $$\mbox{speed-up}=\frac{\mbox{CPU of RK}}{\mbox{CPU of DSBGS}}.$$ In each experiment, ITER and CPU are averaged over 20 trials.

\subsection{Synthetic data}   

Two types of coefficient matrices are generated as follows. 
\bit
\item[$\bullet$] Type I: For given $m$, $n$, $r = rank(\mbf A)$, and $\kappa>1$, we construct a matrix $\mbf A$ by $$\bf A = UDV^\rmt,$$ where $\mbf U\in\mbbr^{m\times r} $ and $\mbf V\in \mbbr^{n\times r}$. Entries of $\mbf U$ and $\mbf V$ are generated from a standard normal distribution, and then, columns are orthonormalized: 
$${\tt [U,\sim] = qr(randn(m,r),0);\qquad  [V,\sim] = qr(randn(n,r),0);}$$ The matrix $\mbf D$ is an $r\times r$ diagonal matrix whose diagonal entries are uniformly distributed numbers in $(1,\kappa)$: $${\tt D = diag(1+(\kappa-1).*rand(r,1));}$$ So the condition number of $\mbf A$, which is defined as $\sigma_1(\mbf A)/\sigma_r(\mbf A)$, is upper bounded by $\kappa$.
\item[$\bullet$] Type II: For given $m$, $n$, entries of $\mbf A$  are generated from a standard normal distribution: $$\tt A=randn(m,n);$$ So $\mbf A$ is a full (column or row) rank matrix with probability one.     
\eit 

\begin{figure}
\centerline{\epsfig{figure=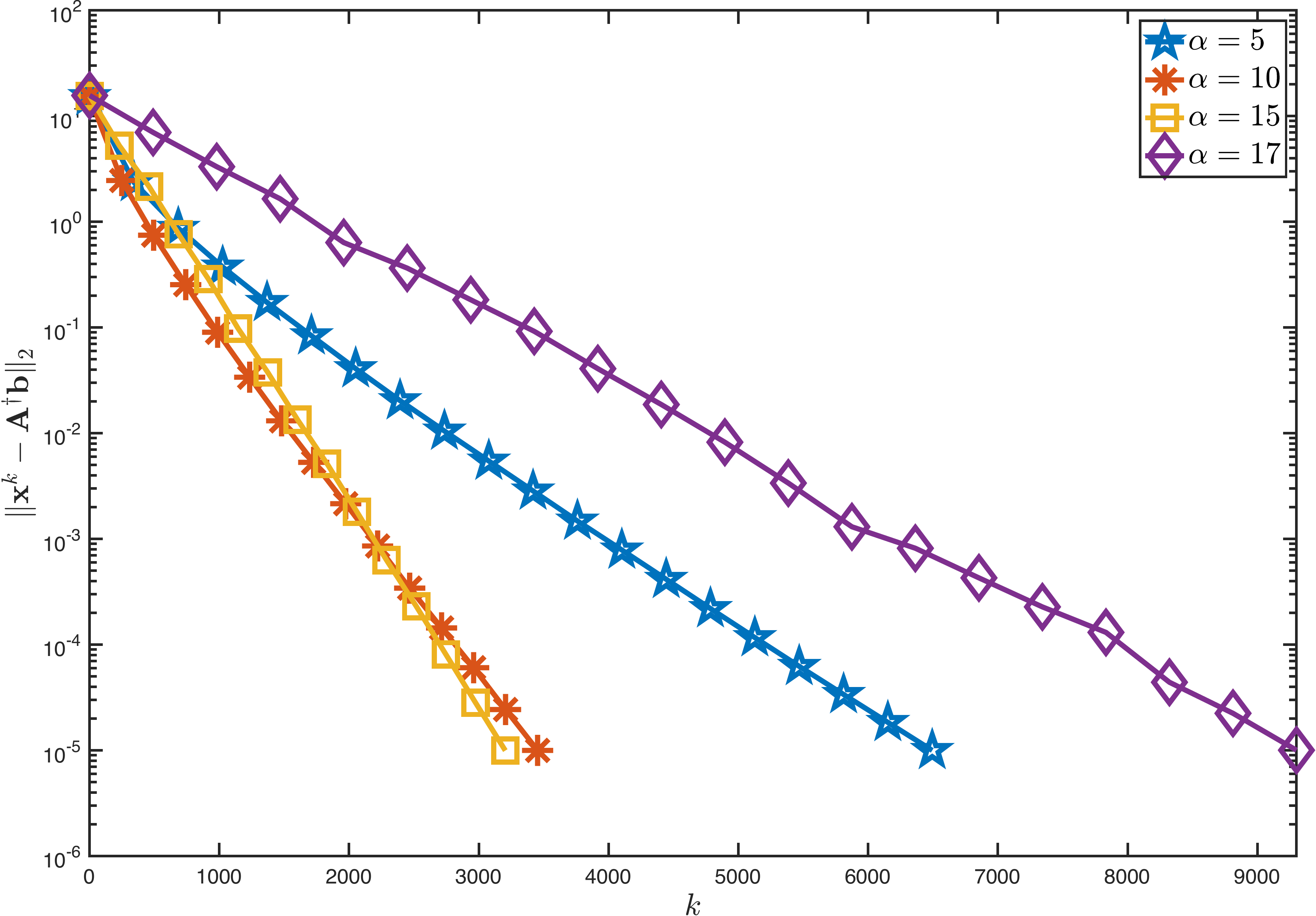,height=3in}}
\caption{The average (20 trials of each case) convergence history of DSBGS with different step size ($\alpha=5,10,15,17$) and fixed block size ($\ell=50,\tau=50$) for a full column rank consistent linear system with random coefficient matrix $\mbf A$ of Type II ({\tt A=randn(500,250)}).}
\label{fig1}       
\end{figure}

\begin{figure}
\centerline{\epsfig{figure=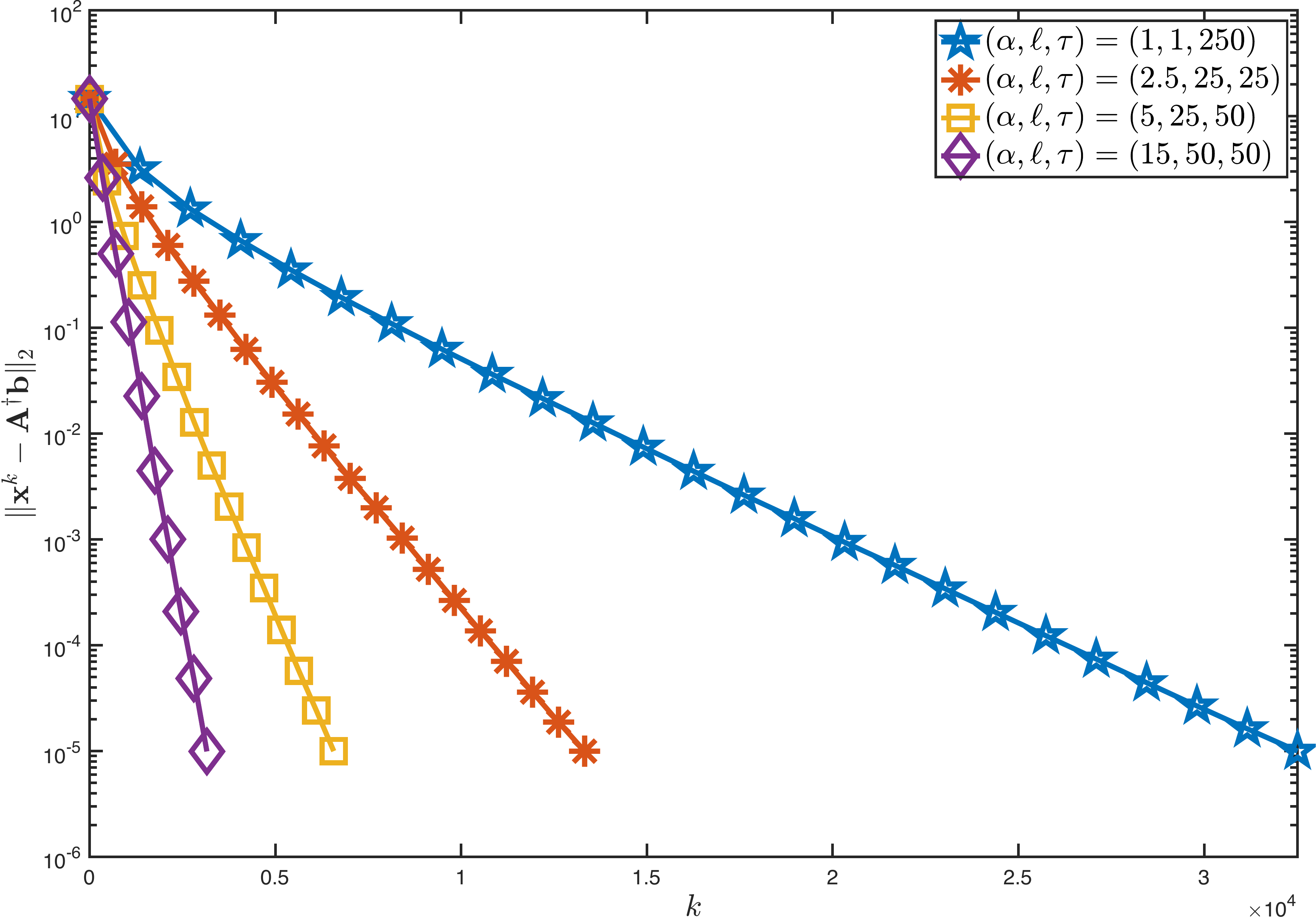,height=3in}}
\caption{The average (20 trials of each case)  convergence history of DSBGS with different step size and different block size  for a full column rank consistent linear system with random coefficient matrix $\mbf A$ of Type II ({\tt A=randn(500,250)}).}
\label{fig2}       
\end{figure}

In Figures \ref{fig1} and \ref{fig2} we plot the error $\|\mbf x^k-\mbf A^\dag\mbf b\|_2$ of DSBGS with different step size $\alpha$ and different block size for full column rank consistent linear systems. From these figures, we observe that appropriate step size and block size improve the convergence remarkably. In Tables \ref{tab1} and \ref{tab2}, we report the average numbers of iterations and the average  computing times for RK and DSBGS. From these results we see DSBGS vastly outperforms RK in terms of computing times with significant speed-ups for general (overdetermined or underdetermined, full column rank or rank deficient) consistent linear systems. It should be noted that the step size $\alpha\in(0,2/{(t\beta)})$ given in Theorem \ref{main3} is a sufficient condition for DSBGS's convergence. Numerical experiments show that DSBGS with some large step size (for example $\alpha=15$ in the experiment for Figure 2) converges much faster than with step size $\alpha\in(0,2/{(t\beta)})$. 

\begin{table}
\caption{The average (20 trials of each experiment)  {\rm ITER} and {\rm CPU} of {\rm RK} and {\rm DSBGS}($\alpha,\ell,\tau$) for consistent linear systems  with  random coefficient matrices $\mbf A$ of Type I: ${\bf A=UDV^\rmt}$.}
\label{tab1}       
\begin{center}
\begin{tabular}{ccr|rr|rrr|c} \toprule $m\times n$ & rank & $\kappa$ & \multicolumn{2}{|c|}{RK: ITER,  CPU} & \multicolumn{3}{|c|}{DSBGS: ITER,  CPU,   $(\alpha,\ell,\tau)$} & speed-up\\ \noalign{\smallskip} \hline \noalign{\smallskip}
$125\times 250$ & 100 & 2  & 3162.55 & 0.0805 & 628.85 & 0.0258 & $(5,5,n)$ & 3.12 \\  
$125\times 250$ & 100 & 10 & 26413.55 & 0.6813 & 5255.70 & 0.1751 & $(5,5,n)$ & 3.89  \\  
$125\times 250$ & 100 & 20 & 158506.00 & 3.8662 &31751.35 & 1.0482 & $(5,5,n)$ & 3.69 \\  
$250\times 500$ & 200 & 2  & 6791.10 & 0.2270 & 638.40 & 0.0482 & $(10,10,n)$ & 4.71 \\
$250\times 500$ & 200 & 10 & 69568.35 & 2.3030 & 6912.85 & 0.5141 & $(10,10,n)$ &  4.48 \\  
$250\times 500$ & 200 & 20 & 252768.05 & 8.3663 & 25328.15 & 1.8639 & $(10,10,n)$ & 4.49 \\  
$250\times 125$ & 125 & 2  & 4215.20 & 0.1138 & 975.25 & 0.0291 & $(5,25,25)$ & 3.90 \\  
$250\times 125$ & 125 & 10 & 38675.35 & 1.0213 & 8241.10 & 0.2396 & $(5,25,25)$ & 4.26  \\
$250\times 125$ & 125 & 20 & 101769.10 & 2.6832 & 23687.75 & 0.6682 & $(5,25,25)$ & 4.02 \\  
$500\times 250$ & 250 & 2  & 8637.00 & 0.2758 & 993.80 & 0.0602 & $(10,50,50)$ & 4.58 \\  
$500\times 250$ & 250 & 10 & 85328.80 & 2.7336 & 9732.80 & 0.5871 & $(10,50,50)$ & 4.66 \\  
$500\times 250$ & 250 & 20 & 448211.30 & 14.1407 & 45301.45 & 2.7288 & $(10,50,50)$ & 5.18 \\
\bottomrule
\end{tabular}
\end{center}
\end{table}

\begin{table}
\caption{The average (20 trials of each experiment)  {\rm ITER} and {\rm CPU} of {\rm RK} and {\rm DSBGS}($\alpha,\ell,\tau$) for consistent linear systems  with  random coefficient matrices $\mbf A$ of Type II: {\tt A=randn(m,n)}. Here $\kappa(\mbf A)=\sigma_1(\mbf A)/\sigma_r(\mbf A)$.} 
\label{tab2}       
\begin{center}
\begin{tabular}{rc|rr|rrr|c} \toprule $m\times n$  & $\kappa(\mbf A)$ & \multicolumn{2}{|c|}{RK: ITER, CPU} & \multicolumn{3}{|c|}{DSBGS: ITER, CPU, $(\alpha,\ell,\tau)$} & speed-up \\ \noalign{\smallskip}  \hline \noalign{\smallskip}
$125\times 250$   & 5.66  & 16118.20 & 0.4160 & 3210.75 & 0.1063 & $(5,5,n)$ & 3.91 \\  
$125\times 500$   & 2.87 & 5876.40 & 0.1719 & 1134.60  & 0.0756 & $(5,5,n)$ & 2.28 \\  
$125\times 1000$  & 2.09 & 3867.00 & 0.1571 & 727.65 & 0.0691 & $(5,5,n)$ & 2.27 \\  
$250\times 500$   & 5.56 & 30557.05 & 1.0139 & 3091.35 & 0.2273 & $(10,10,n)$ & 4.46 \\
$250\times 1000$  & 2.96 & 12015.30 & 0.5085 & 1174.55 & 0.1310 & $(10,10,n)$ & 3.88 \\
$250\times 2000$  & 2.06  & 7927.15 & 0.4963 & 751.75 & 0.1959 & $(10,10,n)$ & 2.53 \\  
$500\times 750$   & 9.66  & 173700.40 & 7.2715 & 17381.35 & 1.7435 & $(10,10,n)$ & 4.17 \\
$500\times 1500$  & 3.66  & 33019.55 & 2.1238 & 3325.55 & 0.5445 & $(10,10,n)$ & 3.90 \\
$500\times 3000$  & 2.34 & 18520.60 & 2.4642 & 1813.70 & 0.7393 & $(10,10,n)$ & 3.33 \\  
$250\times 125$   & 5.33 & 13981.95 & 0.3709 & 3053.15 & 0.0864 & $(5,25,25)$ & 4.29 \\  
$500\times 125$   & 2.96 & 6095.05 & 0.1811 & 1338.45 & 0.0393 & $(5,25,25)$ & 4.60 \\  
$1000\times 125$  & 2.07 & 4112.30 & 0.1488 & 1003.05 & 0.0344 & $(5,25,25)$ & 4.32 \\  
$500\times 250$   & 5.77 & 32563.80 & 1.0330 & 6958.60 & 0.4097 & $(5,50,25)$ & 2.52 \\
$1000\times 250$  & 2.87 & 11965.45 & 0.4537 & 2724.85 & 0.1763 & $(5,50,25)$ & 2.57 \\
$2000\times 250$  & 2.12 & 8489.30 & 0.4305 & 2086.25 & 0.1430 & $(5,50,25)$ & 3.01 \\  
$750\times 500$   & 9.50 & 148619.65 & 5.9663 &  32345.90 & 2.8518 & $(5,50,50)$ & 2.09  \\
$1500\times 500$  & 3.68 & 33655.05 & 1.6979 & 7212.25 & 0.6977 & $(5,50,50)$ & 2.43 \\
$3000\times 500$  & 2.36 & 18990.70 & 1.3557 & 4378.90 & 0.5279 & $(5,50,50)$ & 2.57 \\  
\bottomrule
\end{tabular}
\end{center}
\end{table}

\subsection{Real-world data}

We test RK and DSBGS on eight real-world problems from the University of Florida sparse matrix collection \cite{davis2011unive}: {\tt abtaha1}, {\tt WorldCities}, {\tt cari}, {\tt df2177}, {\tt flower\_5\_1}, {\tt football}, {\tt relat6}, {\tt Sandi\_authors}. The first two matrices are of full column rank and the last six matrices are rank-deficient. In Table \ref{tab3}, we report the average numbers of iterations and the average computing times of RK and DSBGS. We observe that DSBGS based on good choices of step size and block size significantly outperforms RK. Moreover, good step size and block size are problem dependent.

\begin{table}
\caption{The average (20 trials of each experiment) ITER and {\rm CPU}  of  {\rm RK} and {\rm DSBGS}($\alpha,\ell,\tau$) for consistent linear systems  with coefficient matrices from the University of Florida sparse matrix collection. Here $\kappa(\mbf A)=\sigma_1(\mbf A)/\sigma_r(\mbf A)$. The first two matrices are of full column rank and the last six matrices are rank deficient.}
\label{tab3}       
\begin{center}
{\footnotesize\begin{tabular}{lrr|rr|rrr|c} \toprule Matrix & $m \times n$  & $\kappa(\mbf A)$ & \multicolumn{2}{|c|}{RK: ITER, CPU} & \multicolumn{3}{|c|}{DSBGS: ITER, CPU, $(\alpha,\ell,\tau)$}  & speed-up
\\ \noalign{\smallskip} \hline \noalign{\smallskip}
{\tt abtaha1} & $14596\times209$  & 12.23 & 1.83e05 & 42.2491 & 3.91e04 & 2.8988 & $(5,10,n)$ & 14.57 \\
{\tt WorldCities}  & $315\times 100$  & 66.00& 7.38e04 & 1.9794 & 3.09e04 & 0.8624 & $(2.5,10,n)$ & 2.30 \\  
{\tt cari} & $400\times1200$  & 3.13 & 9.79e03 & 0.5035 & 2.76e03 & 0.3186 & $(2.5,5,n)$ & 1.58 \\  
{\tt df2177} & $630\times10358$ & 2.01 & 1.63e04 & 8.5768 & 2.66e03 & 5.0829 & $(5,10,n)$ & 1.69 \\  
{\tt flower\_5\_1}  & $211\times201 $  & 13.70 & 9.55e04 & 2.6053 & 3.83e04 & 1.2001 & $(2.5,5,n)$ & 2.17 \\ 
{\tt football}  & $35\times35 $ & 166.47 & 7.88e05 & 15.1897 & 3.94e05 & 8.4485 & $(2,4,n)$ & 1.80 \\ 
{\tt relat6}  & $2340\times 157$  & 7.74 & 2.66e04 & 1.4331 & 1.03e04 & 0.3969 & $(2.5,10,n)$ & 3.61 \\ 
{\tt Sandi\_authors}  & $86\times86 $   & 189.58 &  2.16e06 & 45.6432 & 8.66e05 & 21.5159 & $(2.5,5,n)$ & 2.12 \\  
   \bottomrule
\end{tabular}}
\end{center}
\end{table}

\section{Concluding remarks} 
We have proposed a doubly stochastic block Gauss--Seidel algorithm for solving linear systems and prove its convergence theory. The randomized Kaczmarz algorithm, the randomized Gauss--Seidel algorithm, and the doubly stochastic Gauss--Seidel algorithm are special cases of the doubly stochastic block Gauss--Seidel algorithm. Numerical experiments show that appropriate step size and block size significantly improve the performance. Finding appropriate variable step size, proposing more effective sampling strategies for submatrices, and designing other block variants via the ideas in \cite{necoara2019faste} should be valuable topics in the future study.

\section*{Acknowledgments}
The research of the first author was supported by the National Natural Science Foundation of China (No.11771364) and the Fundamental Research Funds for the Central Universities (No.20720180008).


\end{document}